\documentclass[11pt]{article}
\usepackage{amsthm,amsfonts,amsmath,amssymb}
\newtheorem{thm}{Theorem}[section]

\newtheorem{lem}[thm]{Lemma}
\newtheorem{cor}[thm]{Corollary}

\def\di{\bigm|} \def\lg{\langle} \def\rg{\rangle}

\def\f{\noindent}
\def\maxsgp{<\!\!\!\cdot\ }

\def\mod{\hbox{\rm mod }}

\def\demo{\f{\bf Proof}\hskip10pt}

 \def\og{\overline G} \def\oh{\overline H}
  
\def\ok{\overline K}  
\def\on{\overline N}

\def\qed{\hfill $\Box$}

\def\M{{\rm M}}

\def\lg{\langle}
\def\rg{\rangle}

\def\rr#1{\item[{\rm (#1)}]}

\def\A{\mathcal{A}$$}

\begin{document}
\title{On Finite Metahamiltonian $p$-Groups
\thanks{This work was supported by NSFC (no. 11371232),
by NSF of Shanxi Province (no. 2013011001-1) and Shanxi Scholarship
Council of China (No. [2011]8-059). }}
\author{
Lijian An, Qinhai Zhang\thanks{Corresponding author. e-mail: zhangqh@dns.sxnu.edu.cn}\\
Department of Mathematics, Shanxi Normal University\\
Linfen, Shanxi, 041004 PR China}

\date{}
\maketitle

\begin{abstract}
A group is called metahamiltonian if all non-abelian subgroups of it
are normal. This concept is a natural generalization of Hamiltonian
groups. In this paper, the properties of finite metahamiltonian
$p$-groups are investigated.

\medskip
\noindent{\bf Keywords} Dedekindian groups, metahamiltonian groups,
$\mathcal{A}_2$-groups

\noindent{\it 2000 Mathematics subject classification:} 20D15.
\end{abstract}

\baselineskip=16pt

\section{Introduction}

A group is called Dedekindian if every subgroup of it is normal. In
1897, Dedekind classified finite Dedekindian groups in \cite{Ded}. In
1933, Baer classified infinite Dedekindian groups in \cite{Baer}. A
non-abelian Dedekindian group is also called Hamiltonian.

A non-abelian group is called metahamiltonian if all non-abelian
subgroups of it are normal. This concept is a natural generalization of
Hamiltonian groups. In the 1960's and 70's, many scholars researched
metahamilton groups. Romalis and Sesekin \cite{RS1,RS2,RS3}
investigated some properties on infinite metahamiltonian groups, and
Nagrebeckii\cite{Nag1, Nag2, Nag3} studied finite metahamiltonian
groups. Nagrebeckii \cite{Nag2} proved the following theorem:

\begin{thm}\label{non-nilpotent}
Suppose that $G$ is a finite non-nilpotent group. Then $G$ is
metahamiltonian if and only if $G=SZ(G)$ where $S$ is one of the
following groups:

{\rm (1)} $P\rtimes Q$, where $P$ is an elementary $p$-group, $Q$ is
cyclic and $(p,|Q|)=1$;

{\rm (2)} $Q_8\rtimes Q$, where $Q$ is cyclic and $(|Q|,2)=1$;

{\rm (3)} $P\rtimes Q$, where $|P|=p^3, p\geq 5$, $Q$ is cyclic and
$(p,|Q|)=1$.
\end{thm}
In \cite{Nag2}, more detailed information on $S$ is given. Since a
nilpotent group is the direct product of its Sylow subgroups, by the
above theorem, to study finite metahamiltonian groups, we only need
consider finite metahamiltonian $p$-groups, which is more complex
than the situation of non-nilpotent.

Metahamiltonian $p$-groups contain many important classes of
$p$-groups. For example, finite $p$-groups all of whose subgroups of
index $p^2$ are abelian, are metahamiltonian. All such groups are
determined, see \cite{BJ,Ber2,Dra,Kaz,She,ZSAX} for the
classification. Another example is finite $p$-groups all of whose
non-normal subgroups are cyclic. See \cite{Pas}. The study of metahamilton $p$-groups is an old problem and many scholar consider it important. 
In this paper, the properties of finite metahamiltonian $p$-groups
are investigated. These properties are useful in the classification of metahamilton $p$-groups
\cite{AF}.

\section{Preliminaries}

Let $G$ be a finite group. $G$ is said to be {\it minimal
non-abelian}, if $G$ is non-abelian, but every proper subgroup of
$G$ is abelian. A finite $p$-group $G$ is called an {\it
$\A_t$-group} if every subgroup of index $p^t$ of $G$ is abelian,
but there is at least one non-abelian subgroup of index $p^{t-1}$.
So ${\cal A}_1$-groups are just the minimal non-abelian $p$-groups.

\medskip

Let $G$ be a finite $p$-group. We define $\Lambda_1(G)=\{ a\in G\di
a^p=1 \}$, $V_1(G)=\{a^{p}\di a\in G \}$, $\Omega_1(G)=\lg
\Lambda_1(G) \rg=\lg a\in G\di a^{p}=1\rg$, and $\mho_1(G)=\lg
V_1(G)\rg=\lg a^{p}\di a\in G\rg$; $G$ is called {\it $p$-abelian}
if $(ab)^{p}=a^{p}b^{p}$ for all $a,b\in G$. We use $c(G)$ and
$d(G)$ to denote the nilpotency class and minimal number of
generators, respectively.

We use $M_p(m,n)$ to denote groups $\langle a, b \di
a^{p^m}=b^{p^n}=1, a^b=a^{1+p^{m-1}}\rangle$, where $m\geq 2$, and
use $M_p(m,n,1)$ to denote groups $\langle a, b, c \di
a^{p^m}=b^{p^n}=c^p=1, [a, b]=c, [c, a]=[c, b]=1 \rangle$, where
 $m+n\geq 3$ for $p=2$ and $m\geq n$. We can give a presentation of
minimal non-abelian $p$-groups as follows:

\begin{thm}\label{thm=Redei}{\rm (\cite {R})}{\rm(R\'edei)}
Let $G$ be a minimal non-abelian $p$-group. Then $G$ is $Q_8$,
$M_p(m,n)$, or $M_p(m,n,1)$.
\end{thm}

\medskip

We use $C_n$ and  $C_n^m$  to denote the cyclic group and the direct
product of $m$ cyclic groups of order $n$, respectively; and use
$H*K$ to denote a central product of $H$ and $K$. For undefined
notation and terminology the reader is referred to \cite{Hup}.

\medskip

We have the following information about minimal non-abelian
$p$-groups.

\begin{thm}\label{thm=inner abel properties}{\rm (\cite[Lemma 2.2] {Alj})}
Let $G$ be a finite $p$-group. Then the following conditions are
equivalent:

{\rm (1)} $G$ is an inner abelian $p$-group;

{\rm (2)} $d(G)=2$ and $|G'|=p$;

{\rm (3)} $d(G)=2$ and $Z(G)=\Phi (G)$.

\end{thm}

\begin{lem}{\rm (\cite[$p_{136}$, Proposition 10.28]{Ber1})}
\label{generated by A_1 subgroups} A non-abelian
 $p$-group is generated by minimal non-abelian subgroups.
\end{lem}

Many scholars studied and classified $\mathcal{A}_2$-groups, see,
for example \cite{BJ,Ber2,Dra,Kaz,She,ZSAX}. We have following
Lemma.

\begin{lem}\label{A_2}{\rm(\cite{ZSAX})} Suppose that $G$ is an $\mathcal{A}_2$-group. Then $G$ is
one of the following groups:
\begin{enumerate}
\rr{I} $d(G)=2$ and $G$ has an abelian maximal subgroup.
\begin{enumerate}
\rr{1} $\lg a,b\di a^{8}=b^{2^m}=1, a^b=a^{-1}\rg$, where $m\ge 1$;

\rr{2} $\lg a,b\di a^{8}=b^{2^m}=1, a^b=a^{3}\rg$, where $m\ge 1$;

\rr{3} $\lg a,b\di a^{8}=1, b^{2^m}=a^{4}, a^b=a^{-1}\rg$, where
$m\ge 1$;

\rr{4} $\lg a_1,b\di
a_1^p=a_2^p=a_3^p=b^{p^m}=1,[a_1,b]=a_{2},[a_2,b]=a_3,[a_3,b]=1,[a_i,a_j]=1\rg,$
where $p\ge 5$ for $m=1$, $p\ge 3$ and $1\leq i,j\leq 3$;

\rr{5} $\lg a_1,b\di
a_1^p=a_2^p=b^{p^{m+1}}=1,[a_1,b]=a_{2},[a_2,b]=b^{p^m},[a_1,a_2]=1\rg$,
where $p\ge 3$;

\rr{6} $\lg a_1,b\di
a_1^{p^2}=a_2^p=b^{p^m}=1,[a_1,b]=a_{2},[a_2,b]=a_1^{\nu
p},[a_1,a_2]=1\rg,$ where $p\ge 3$ and $\nu=1$ or a fixed quadratic
non-residue modulo $p$.

\rr{7} $\lg a_1,a_2,b\di
a_1^9=a_2^3=1,b^3=a_1^3,[a_1,b]=a_2,[a_2,b]=a_1^{-3},[a_2,a_1]=1\rg$.

\end{enumerate}

\rr{II} $d(G)=3$, $|G'|=p$ and $G$ has an abelian maximal subgroup.

\begin{enumerate}

\rr{8} $\langle {a},{b},{x} \di {a}^4={x}^{2}=1,
{b}^2={a}^2=[{a},{b}], [{x},{a}]=[{x},{b}]=1 \rangle \cong Q_8
\times C_{2}$;

\rr{9}  $\langle {a}, {b},{x} \di {a}^{p^{n+1}}={b}^{p^m}={x}^{p}=1,
[{a},{b}]={a}^{p^{n}}, [{x},{a}]=[{x},{b}]=1\rangle\cong
M_p(n+1,m)\times C_{p}$;

\rr{10}  $\langle {a}, {b}, {c},{x} \di
{a}^{p^n}={b}^{p^m}={c}^{p}={x}^p=1, [{a},{b}]={c},
[{c},{a}]=[{c},{b}]=[{x},{a}]=[{x},{b}]=1\rangle\cong
M_p(n,m,1)\times C_{p}$, where $n\ge m$, and $n\ge 2$ if $p=2$;

\rr{11}  $\langle {a}, {b}, {x} \di {a}^4=1,
{b}^2={x}^{2}={a}^2=[{a},{b}], [{x},{a}]=[{x},{b}]=1\rangle \cong
Q_8 \ast C_{4}$;

\rr{12}  $\langle {a}, {b}, {x} \di
{a}^{p^n}={b}^{p^m}={x}^{p^{2}}=1, [{a},{b}]={x}^{p},
[{x},{a}]=[{x},{b}]=1 \rangle\cong M_p(n,m,1)\ast C_{p^{2}}$, where
$n\ge 2$ if $p=2$ and $n\ge m$.

\end{enumerate}

\rr{III} $d(G)=3$, $|G'|=p^2$ and $G$ has an abelian maximal
subgroup.

\begin{enumerate}

\rr{13} $\lg {a},{b},{c}\di {a}^{4}={b}^{4}=1,
{c}^2={a}^2{b}^2,[{a},{b}]={b}^2,[{c},{a}]={a}^2,[{c},{b}]=1\rg$;

\rr{14} $\lg {a},{b},{d}\di {a}^{p^{m}}={b}^{p^2}={d}^{p}=1,
[{a},{b}]={a}^{p^{m-1}},[{d},{a}]={b}^p,[{d},{b}]=1\rg$, where
$m\geq 3$ if $p=2$;

\rr{15} $\lg {a},{b},{d}\di {a}^{p^m}={b}^{p^2}={d}^{p^2}=1,
[{a},{b}]={d}^p,[{d},{a}]={b}^{jp},[{d},{b}]=1\rg$, where $(j,p)=1$,
$p>2$, $j$ is a fixed quadratic non-residue modulo $p$, and $-4j$ is
a quadratic non-residue modulo $p$;

\rr{16} $\lg {a},{b},{d}\di
 {a}^{p^m}={b}^{p^2}={d}^{p^2}=1,[{a},{b}]={d}^p,[{d},{a}]={b}^{jp}{d}^p,[{d},{b}]=1\rg$,
where if $p$ is odd, then $4j =1-\rho^{2r+1}$ with $1\le
r\le\frac{p-1}{2}$ and $\rho$ the smallest positive integer which is
a primitive root $(\mod p)$; if $p = 2$, then $j = 1$.

\end{enumerate}

\rr{IV} $d(G)=2$ and $G$ has no abelian maximal subgroup.

\begin{enumerate}

\rr{17} $\lg a, b\di a^{p^{r+2}}=1,b^{p^{r+s+t}}=a^{p^{r+s}},[a,
b]=a^{p^r}\rg$, where $r\ge 2$ for $p=2$, $r\ge 1$ for $p\ge 3$,
$t\ge 0$, $0\le s\le 2$ and $r+s\ge 2$;

\rr{18} $\lg a,b\di a^{p^2}=b^{p^2}=c^p=1,[a,b]=c,[c,a]=b^{\nu p},
[c,b]=a^{p}\rg$, where $p\ge 5$, $\nu$ is a fixed square non-residue
modulo $p$;

\rr{19} $\lg a,b\di a^{p^2}=b^{p^2}=c^p=1,[a,
b]=c,[c,a]=a^{-p}b^{-lp},[c, b]=a^{-p}\rg$, where $p\ge 5$,
$4l=\rho^{2r+1}-1$, $r=1, 2, \dots, \frac{1}{2}(p-1)$, $\rho$ is the
smallest positive integer which is a primitive root modulo $p$;

\rr{20} $\lg a, b\di a^9=b^9=c^3=1,[a,b]=c,[c,a]=b^{-3}, [c,
b]=a^3\rg$;

\rr{21} $\lg a, b\di a^9=b^9=c^3=1, [a,b]=c, [c,a]=b^{-3}, [c,
b]=a^{-3}\rg$.

\end{enumerate}

\rr{V} $d(G)=3$ and $G$ has no abelian maximal subgroup.

\begin{enumerate}

\rr{22} $\lg a,b,d\di a^4=b^4=d^4=1, [a,b]=d^2, [d,a]=b^2d^2,
[d,b]=a^2b^2,[a^2,b]=[b^2,a]=1\rg$.

\end{enumerate}
\end{enumerate}

\end{lem}

Analyzing the group list in Lemma {\rm\ref{A_2}}, we have following
lemma.

\begin{lem}\label{A_2-property} Suppose that $G$ is an $\mathcal{A}_2$-group with order $p^n$.

{\rm (1)} $d(G)\le 3$ and $c(G)\le 3$;

{\rm (2)} If $d(G)=2$ and $\exp(G')=p$, then $c(G)=3$.

{\rm (3)} If $c(G)>2$ and $\exp(G')=p$, then $d(G)=2$ and $p$ is
odd.

\end{lem}

\begin{thm}{\rm (\cite[Statz 6.5]{Hup}}\label{thm=engle}
If $[x,y,y]=1$ for all $x,y\in G$, then $G$ is nilpotent and
$c(G)\le 3$. In addition, if $G$ has no element of order $3$, then
$c(G)\le 2$.
 \end{thm}
A finite $p$-group $G$ is called {\em metacyclic} if it has a cyclic
normal subgroup $N$ such that $G/N$ is also cyclic. 

\begin{lem}{\rm (\cite{Bla})}
\label{metacyclic} Suppose that $G$ is a finite $p$-group. Then $G$
is metacyclic if and only if $G/\Phi(G')G_3$ is metacyclic.
\end{lem}

\section{Properties of finite
metahamiltonian $p$-groups}

\begin{thm}\label{herit}
Let $G$ be a finite metahamiltonian $p$-group. Then sections of $G$
are all metahamiltonian.
\end{thm}
\demo It is straight forward. \qed

\begin{thm}\label{thm=pandun}
Let $G$ be a finite $p$-group. Then $G$ is metahamiltonian if and
only if every minimal non-abelian subgroup is normal in $G$.
\end{thm}
\demo If $G$ is metahamiltonian, then, by the definition of
metahamiltonian, every minimal non-abelian subgroup is normal in
$G$. On the other hand, if every minimal non-abelian subgroup is
normal in $G$, then, by Lemma \ref{generated by A_1 subgroups},
every non-abelian subgroup is normal in $G$.\qed

\begin{thm}\label{Closure of x}
Let $G$ be a finite metahamiltonian $p$-group. Then, for all $x\in
G$, $\lg x\rg^G$ is abelian or minimal non-abelian.
\end{thm}
\demo Suppose that $\lg x\rg^G$ is not abelian. Then there exists
$g\in G$ such that $[x,x^g]\ne 1$. Let $K=\lg x,x^g\rg$. Then $K$ is
normal in $G$ since $G$ is metahamiltonian. Hence $K=\lg x\rg^G$.
Let $y=x^g$ and $L=\lg x,x^y\rg=\lg x,[x,y]\rg$. Then $L<K$ and
hence $L$ is not normal in $G$. It follows that $L$ is abelian. That
is, $[x,y,x]=1$. Since $\lg y\rg^G=\lg x^g\rg^G=\lg x\rg^G$,
similarly we have $[x,y,y]=1$. Hence $c(K)=2$.

Let $S=\lg x,y^p\rg$. Then $S<K$ and hence $S$ is not normal in $G$.
It follows that $S$ is abelian and hence $[x,y^p]=1$. Since
$c(K)=2$, we get $[x,y]^p=1$. Thus $K'=\lg [x,y]\rg$ is of order
$p$. By Theorem \ref{thm=inner abel properties}, $K$ is minimal
non-abelian. \qed

\begin{thm}\label{thm=metahamilton p-gp G, c(G) is at most 3}
Let $G$ be a metahamiltonian $p$-group. Then $c(G)\le 3$. In
particular, $G$ is metabelian.
\end{thm}

\demo By Theorem \ref{Closure of x}, for all $x\in G$, $K=\lg x
\rg^G$ is abelian or minimal non-abelian. Then $K'=1$ or $|K'|=p$.
Since $K'\unlhd G$, we get $K'\le Z(G)$. Let $\bar{G}=G/Z(G)$. Then,
for all $\bar{x}\in\bar{G}$, $\lg \bar{x}^{\bar{G}}\rg$ is abelian.
Hence $\bar{G}$ satisfies the 2-Engel condition. By Theorem
\ref{thm=engle}, $c(\bar{G})=2$ for $p\ne 3$ and $c(\bar{G})\le 3$
for $p=3$. It follows that $c({G})\le 3$ for $p\ne 3$ and $c({G})\le
4$ for $p=3$.

We claim that $c(G)\le 3$. If not, then $p=3$ by the above argument. Let
$G$ be a counterexample with minimal order. By Theorem \ref{herit},
$c(G)=4$, $|G_4|=p$ and the nilpotency class of every proper section
of $G$ is at most 3. Hence we may assume that $G_4=\lg
[a,b,c,d]\rg$, where $a,b,c,d\in G\setminus\Phi(G)$. Let
$x=[a,b,c]$. Then $N=\lg x,d\rg$ is minimal non-abelian by Theorem
\ref{thm=inner abel properties}. By hypothesis, every subgroup which
contains $N$ is normal in $G$. It follows that $G/N$ is Dedekindian.
Since $p=3$, $G/N$ is abelian. It follows that $G'\le N$. Since
$d\not\in\Phi(G)$, we have $N\cap \Phi(G)<N$ and hence $G'\le
N\cap\Phi(G)<N$. It follows that $G'$ is abelian. Then
$[[c,d],[a,b]]=1$. Since $[a,b]\in G'<N$ and $d\in N$, $[d,[a,b]]\in
N'\le Z(G)$. It follows that $[d,[a,b],c]=1$. By Witt's
formula, we have $[[a,b],c,d]=1$, a contradiction.\qed

\medskip

\begin{thm}\label{G'<N}Let $G$ be a finite $p$-group. $G$ is metahamiltonian if and only if
$G'$ is contained in every non-abelian subgroup of $G$.
\end{thm}
\demo If $G'$ is contained in every non-abelian subgroup of $G$,
then every non-abelian subgroup of $G$ is normal in $G$. Hence
sufficiency holds. In the following we prove the necessity.

Let $G$ be a counterexample with minimal order. Then $G$ is
metahamiltonian and there exists a minimal non-abelian subgroup
$N=\lg a,b\rg$ such that $G'\not\le N$. Since $G$ is
metahamiltonian, subgroups containing $N$ are normal in $G$. Hence
$G/N$ is Hamiltonian.

By the minimality of $G$, $G/N\cong Q_8$. Let $G/N=\lg xN,yN\rg$ and
$H=\lg x,y\rg$. Then $G=HN$, $H/(H\cap N)\cong Q_8$,
$z:=[x,y]\not\in N$, $H\cap N\le\Phi(H)$ and $H\cap N=\lg
x^4,x^2y^2,x^2[x,y]\rg^H$. Since $z\in\lg x\rg^H$, it follows from
Theorem \ref{Closure of x} that $\lg z,x\rg$ is abelian or minimal
non-abelian. Hence $[z,x^2]=[z,x]^2=1$. The same reason gives that
$[z,y^2]=[z,y]^2=1$ and hence $\exp(H_3)\le 2$. Since $\Phi(H)=\lg
x^2,y^2,H'\rg$ and $H'$ is abelian (by Theorem \ref{thm=metahamilton
p-gp G, c(G) is at most 3}), we have $[\Phi(H),z]=1$. In particular,
$[H\cap N, z]=1$. In the following, we deduce a contradiction on five
cases:

\medskip

Case 1. $H\cap N=N$.

\medskip

In this case, $[N,z]=1$. Let $M=\lg za,b\rg$. Then Theorem
\ref{thm=inner abel properties} gives that $M$ is minimal
non-abelian, and hence $G/M$ is also Dedekindian. Since $z\not\in
M$, $G/M$ is not abelian. By the minimality of $G$, $H/M=G/M\cong
Q_8$. It follows that $M=\lg x^4,x^2y^2,x^2[x,y]\rg^H=N=\lg a,b\rg$,
a contradiction.

\medskip

Case 2. $H\cap N<N$ and $H\cap N\not\le\Phi(N)$.

\medskip

In this cases, $H\cap N$ contains a generator of $N$. Without
losing generality, we assume that $a\in H\cap N$ and $b\not\in H\cap
N$. Then $[z,a]=1$. Since $H\cap N$ is abelian, we have
$[x^2y^2,x^2[x,y]]=1$, and hence $[x^2,y^2]=1$. By calculation, we
have $[x^2,y^2]=[x^2,y]^2=[x,y]^4=z^4$. If $z^2\neq 1$,  then $\lg
z^2\rg=\mho_1(H')$ is a minimal normal subgroup of $G$. Hence we
have $z^2\in Z(G)$. Particularly, $[z,b]^2=[z^2,b]=1$.

Subcase 2.1. $[z,b]\neq [a,b]$.

Let $M=\lg za,b\rg$. By Theorem \ref{thm=inner abel properties}, $M$
is minimal non-abelian and hence $G/M$ is also Dedekindian. Since
$z\not\in M$, $G/M$ is not abelian. By the minimality of $G$, we
have $G/M=HM\cong H/(H\cap M)\cong Q_8$. It follows that $H\cap
M=\lg x^4,x^2y^2,x^2[x,y]\rg^H=H\cap N$, and hence $a\in H\cap
N=H\cap M\le M$. Thus $z=(za)a^{-1}\in M$, a contradiction.

Subcase 2.2. $[z,b]=[a,b]$.

Let $L=\lg z,b\rg\cap N$. Then $L$ is normal in $G$. Let $K$ be a
maximal subgroups of $N$ which contains $L$ such that $K\unlhd G$.
Then $G/K$ is of order $2^4$, has two generators, and has a quotient
group which is isomorphic to $Q_8$. By the classification of groups
of order $2^4$, $G/K=\lg xK,yK\rg:=\lg \bar{x},\bar{y}\rg\cong
M_2(2,2)$, which has definition relations $\bar{x}^4=\bar{y}^4=1$
and $[\bar{x},\bar{y}]=\bar{x}^2$. Obviously, $\lg\bar{y}\rg$ and
$\lg\bar{x}\bar{y}\rg$ are not normal in $G/K$. It follows that
their complete inverse images are also not normal in $G$, hence are
abelian. It follows that $[y,K]=1$, $[xy,K]=1$. Thus $[H,K]=1$,
which is contrary to $[z,b]=[a,b]\neq 1$.

\medskip

Case 3. $H\cap N<\Phi(N)$.

\medskip

We claim that $H\cap N\neq 1$. Otherwise, $G=H\times N$. Since
$N\cong G/H$ is Dedekind, we have $N\cong Q_8$. In this case, $\lg
xa,yb\rg\cong Q_8$ is not normal in $G$, a contradiction.

We  claim that $N'\le H\cap N$. Otherwise, $G/(H\cap N)$ is also a
counterexample, which is contrary to the minimality of $G$.

Let $\og=G/(H\cap N)$, $\oh=H/(H\cap N)=\lg \bar{x},\bar{y}\rg$ and
$\on=N/(H\cap N)=\lg\bar{a}\rg\times\lg\bar{b}\rg$. Then
$\og=\oh\times \on$ and $\exp(\on)\ge 4$. Without loss of
generality, we may assume that $o(\bar{a})\ge 4$. Let $\ok=\lg
\bar{x}\bar{a}\rg\times\lg\bar{b}\rg$. Then $\ok$ is not normal in
$\og$. It follows that its complete inverse image is not normal in
$G$. Hence $[xa,b]=1$. That is, $[x,b]=[a,b]$. The same reason gives
that $[y,b]=[a,b]$ and $[xy,b]=[a,b]$, a contradiction.

\medskip

Case 4. $H\cap N=\Phi(N)=N'$.

\medskip

In this case, $|N|=2^3$, $|H|=2^4$, $|G|=2^6$ and
$G/N'=H/N'\times\lg aN'\rg\times\lg bN'\rg$. Since $\lg aN'\rg$ and
$\lg bN'\rg$ are normal in $G/N'$, $A:=\lg a,N'\rg$ and $B:=\lg
b,N'\rg$, their complete inverse images, are also normal in $G$.
Noting that $A$ and $B$ are of order $4$, the NC-Theorem gives that
$C_G(A)$ and $C_G(B)$ are maximal in $G$. Let $K=C_G(A)\cap C_G(B)$.
Then $|K|\ge 2^4$. Since $K\cap N=Z(N)=N'$, we have
$|KN|=(|K||N|)/|K\cap N|\ge 2^6$, and hence $G=K*N$. Since
$KN/N\cong K/K\cap N\cong Q_8$, without loss of generality we may
assume that $H=K$. By the classification of groups of order $2^4$,
$H=\lg x,y\rg\cong M_2(2,2)$, which has definition relations
$x^4={y}^4=1,[{x},{y}]={x}^2$ and $N'=H\cap N=\lg x^2y^2\rg$.

Without loss of generality, we may assume that $a\in N$ is of order
$4$. Then $a^2=x^2y^2$. By calculations, we have $[x,ay]=x^2$ and
$(ay)^2=x^2$. It follows that $\lg x,ay\rg$ is neither abelian nor
normal in $G$, a contradiction.

\medskip

Case 5. $H\cap N=\Phi(N)\neq N'$.

\medskip

 Let $\og=G/K$, $\oh=H/K=\lg\bar{x},\bar{y}\rg$ and
$\on=N/K=\lg \bar{a}\rg\times\lg\bar{b}\rg$, where $K$ is a maximal
subgroup of $H\cap N$ such that $K\unlhd G$. By the minimality of
$G$, $\og'$ is contained in every minimal non-abelian subgroup of
$\og$. Since $G'\not\le N$, we have $\og'\not\le \on$and hence $\on$
is abelian. Without loss of generality, we may assume that
$o(\bar{a})=4$. By the classification of groups of order $2^4$,
$\oh\cong M_2(2,2)$, which has definitions $\bar{x}^4=\bar{y}^4=1$,
$[\bar{x},\bar{y}]=\bar{x}^2$, $\bar{a}^2=\bar{x}^2\bar{y}^2$, and
$\Phi(N)/K=(H\cap N)/K=\lg\bar{x}^2\bar{y}^2\rg$. If $\bar{a}\in
Z(\og)$, then $\lg \bar{x},\bar{a}\bar{y}\rg$ is neither abelian nor
normal in $\og$, a contradiction. Hence $\bar{a}\not\in Z(\og)$. If
$[\bar{a},\bar{x}]=\bar{1}$, then
$[\bar{a},\bar{y}]=\bar{x}^2\bar{y}^2$ and hence $\lg
\bar{a}\bar{x},\bar{y}\rg$ is neither abelian nor normal in $\og$, a
contradiction. Hence $[\bar{a},\bar{x}]=\bar{x}^2\bar{y}^2$. The
same reason gives that
$[\bar{a}\bar{b},\bar{x}]=\bar{x}^2\bar{y}^2$. It follows that
$[\bar{b},\bar{x}]=1$. If $[\bar{b},\bar{y}]\neq \bar{1}$, then
$[\bar{b},\bar{y}]=\bar{x}^2\bar{y}^2$. By calculation, $\lg
\bar{x},\bar{b}\bar{y}\rg$ is neither abelian nor normal in $\og$, a
contradiction. Hence $[\bar{b},\bar{y}]=\bar{1}$.

In this case, it is easy to see that $\lg \bar{x},\bar{b}\rg$ and
$\lg \bar{a}\bar{x},\bar{b}\rg$ are not normal in $\og$. It follows
that their complete converse images are not normal in $G$, and hence
they are abelian. Thus $[x,b]=1$ and $[ax,b]=1$, which is contrary
to $[a,b]\neq 1$.\qed

\begin{thm} \label{thm=d(G)=2--metacyclic}Suppose that $G$ is a finite metahamilton $p$-group. If $d(G)=2$ and $\exp(G')>p$,
then $G$ is metacyclic.
\end{thm}
\demo Assume that $G=\lg a,b\rg$ is a counterexample with minimal
order. By Lemma \ref{metacyclic}, $\og:=G/\Phi(G')G_3$ is not
metacyclic. Since $|{\og}'|=p$, ${\og}$ is minimal non-abelian. By
Theorem \ref{thm=Redei}, ${\og}\cong \M_p(n,m,1)$. That is, we may
assume that $\og=\lg \bar{a},\bar{b}\di
\bar{a}^{p^n}=\bar{b}^{p^m}=\bar{c}^p=1,[\bar{a},\bar{b}]=\bar{c},[\bar{c},\bar{a}]=[\bar{c},\bar{b}]=1\rg$.
Since $\lg\bar{a}^p,\bar{b}\rg$, $\lg\bar{b}^p,\bar{a}\rg$,
$\lg(\bar{a}\bar{b})^p,\bar{a}\rg$ and
$\lg(\bar{a}\bar{b})^p,\bar{b}\rg$ are not normal in $\og$, we have
$\lg a^p,b,\Phi(G')G_3\rg$, $\lg b^p,a,\Phi(G')G_3\rg$, $\lg
(ab^p),a,\Phi(G')G_3\rg$ and $\lg (ab)^p,b,\Phi(G')G_3\rg$ are not
normal in $G$. Hence they are all abelian. Thus we have
$\Phi(G')G_3\le Z(G)$ and
$$[a^p,b]=[b^p,a]=[(ab)^p,a]=[(ab)^p,b]=1\eqno(*)$$

If $p=2$, then $(ab)^2=a^2b^2[a,b]$. By $(*)$, $[a,b]\in Z(G)$.
Hence $G'=\lg [a,b]\rg$. By $(*)$, $[a,b]^2=[a^2,b]=1$, which is
contrary to $\exp(G')>2$.

If $p>2$, then, by calculation, we have $[a,b,a]^p=[a^p,b,a]=1$ and
$[a,b,b]^p=[a^p,b,b]=1$. It follows that $\exp(G_3)\le p$. By $(*)$,
$[a,b]^p=[a^p,b]=1$, which is contrary to $\exp(G')>p$.\qed

\begin{lem}\label{lem=not contain metacyclic}
Suppose that $G$ is a finite metahamiltonian $p$-group which has
elementary abelian derived group. If $G$ is not an
$\mathcal{A}_2$-group, then $\mathcal{A}_2$-subgroups of $G$ have
nilpotency class $2$.
\end{lem}
\demo Assume the contrary. Then there exists $K<G$ such that
$K\in\mathcal{A}_2$, $\exp(K')=p$ and $c(K)\ge 3$. Hence $p>2$ and
$K$ is a group of Type (4)--(7) or (18)--(21) in Lemma \ref{A_2}.

\medskip

Case 1: $K=\lg a_1,b\rg$ is a group of Type (4)--(7) in Lemma
\ref{A_2}.

\medskip

Let $H\le G$ such that $K\maxsgp H$. Since $K/Z(K)$ is minimal
non-abelian, non-metacyclic and of order $p^3$, $H/Z(K)$ is not
metacyclic and of order $p^4$. If $d(H/Z(K))=2$, then, by the
classification of groups of order $p^4$, $H/Z(K)$ is of maximal
class. It follows that $K'Z(K)/Z(K)=H_3Z(K)/Z(K)$. Hence $[a_1,b]\in
H_3Z(K)$ and $[a_1,b,b]\in H_4$. Since $[a_1,b,b]\ne 1$, we have
$H_4\ne 1$ and $c(H)\ge 4$, which is contrary to Theorem
\ref{thm=metahamilton p-gp G, c(G) is at most 3}. If $d(H/Z(K))=3$,
then, by the classification of groups of order $p^4$, there exists
$d\in H$ such that $H/Z(K)=K/Z(K)\times \lg dZ(K)\rg$ or
$H/Z(K)=K/Z(K)\ast\lg dZ(K)\rg$. By calculation, $[d^p,k]=[d,k]^p=1$
for all $k\in K$. It follows that $d^p\in Z(K)$ and
$H/Z(K)=K/Z(K)\times \lg dZ(K)\rg$. Since $a_2=[a_1,b]\not\in \lg
a_1,d\rg$, by Theorem \ref{G'<N}, we have $[a_1,d]=1$. The same reason
gives that $[b,d]=1$.  Hence $d\in Z(H)$. In this case $\lg
a_2d,b\rg$ is neither abelian nor normal, a contradiction.

\medskip

Case 2: $K=\lg a_1,b\rg$ is a group of Type (18)--(21) in lemma
\ref{A_2}.

\medskip

Let $H\le G$ such that $K\maxsgp H$. By Theorem \ref{G'<N},
$H'\le\lg c,a\rg\cap\lg c,b\rg=\lg c,a^p,b^p\rg$. It follows that
$H'=K'$ and $H_3=K_3=\lg a^p,b^p\rg$. By the classification of
groups of order $p^4$, there exists $d\in H\setminus K$ such that
$[a,d]\equiv[b,d]\equiv 1 (\mod K_3)$. By calculation,
$[a,d^p]=[a,d]^p=1$ and $[b,d^p]=[b,d]^p=1$. It follows that $d^p\in
Z(K)=K_3$. Since $c\not\in\lg a,d\rg$, by Theorem \ref{G'<N}, we
have $[a,d]=1$. The same reason gives that $[ac,d]=1$. In this case $\lg
a,cd\rg$ is neither abelian nor normal, a contradiction.\qed

\begin{thm}\label{thm=c=3 exp(G')=p} Suppose that $G$ is a finite
metahamiltonian $p$-group having an elementary abelian derived
group. If $c(G)=3$, then $G$ is an $\mathcal{A}_2$-group.
\end{thm}
\demo
 Assume the contrary
and $G$ is a counterexample with minimal order. Then $c(G)=3$ and
$G\in\mathcal{A}_3$.

We claim that $G$ does not satisfy the $2$-Engel condition. If not,
then, by Theorem \ref{thm=engle}, $G$ is a $3$-group. In this case,
there exist $x,y,z\in G$ such that $[x,y,z]\neq 1$. Since $G$ has
minimal order, we have $G=\lg x,y,z\rg$ and $[x,y,z]^3=[x^3,y,z]=1$.
Since $[x,yz,yz]=1$, by calculation, we get $[x,y,z]=[z,x,y]$.
Similar reasons give that $[x,y,z]=[y,z,x]=[z,x,y]$. Let
$[x,y]=c,[y,z]=a,[z,x]=b,[x,y,z]=[y,z,x]=[z,x,y]=d$. Then $G'=\lg
a,b,c,d\rg$. Since $[b,y]=d\neq 1$, we have $\lg b,y\rg\unlhd G$. It
follows that $c=[x,y]\in\lg b,y\rg$. Since $[c,b]=[c,y]=1$, we may
assume that $c=y^{3t}d^w$. Hence
$d=[c,z]=[y^{3t}d^w,z]=[y^{3t},z]=1$, a contradiction.

Since $G$ does not satisfy the $2$-Engel condition, there exist $x,y\in
G$ such that $[x,y,y]\neq 1$. Since $G$ has minimal order, we have
$G=\lg x,y\rg$, $[x,y,y]^p=1$ and $[x,y,x]^p=1$. Let
$[x,y]=c,[c,y]=b$ and $[c,x]=a$. Then $G_3=\lg b,a\rg$ and $G'=\lg
c,G_3\rg$. If $[c,x]\in\lg b\rg$, then, by suitable replacement, we
may assume that $[c,x]=a=1$. Hence we may assume that $\lg
a\rg\cap\lg b\rg=1$.

The maximal subgroups of $G$ are $M=\lg x^iy,\Phi(G)\rg$ and $K=\lg
x,\Phi(G)\rg$, where $i=0,1,\dots,p$. It is easy to see that
$\Phi(G)=\lg x^p, y^p, c,a,b\rg$ is abelian. Since
$[c,x^iy]=a^ib\neq 1$, by Lemma \ref{A_2-property} (2), we have that
$N=\lg c,x^iy\rg\in\mathcal{A}_1$. By Theorem \ref{G'<N}, $G'\le N$.
Since $[cx^p,x^iy]=ba^{i+i {p\choose 2}}\neq 1$, by Lemma
\ref{A_2-property} (2), we have $\lg cx^p,x^iy\rg\in\mathcal{A}_1$.
By Lemma \ref{Closure of x}, $\lg cx^p,x^iy\rg=\lg x^iy\rg^G=N$. It
follows that $x^p\in N$. Since $(x^iy)^p\equiv x^{ip}y^p \ (\mod
G')$, we have $x^{ip}y^p\in N$ and hence $y^p\in N$. Thus
$\Phi(G)\le N$ and $M=N\in\mathcal{A}_1$.

If $[c,x]=a\neq 1$, then, by Lemma \ref{A_2-property} (2), $\lg
c,x\rg\in \mathcal{A}_1$. By Theorem \ref{G'<N}, $G'\le L$. Since
$[cy^p,x]\neq 1$, by Lemma \ref{A_2-property} (2), $\lg cy^p,x\rg\in
\mathcal{A}_1$. By Lemma \ref{Closure of x}, $\lg cy^p,x\rg=\lg
x\rg^G=\lg c,x\rg$. It follows that $y^p\in \lg c,x\rg$ and hence
$\Phi(G)\le \lg c,x\rg$. Thus $K=\lg c,x\rg\in\mathcal{A}_1$.

If $[c,x]=a=1$ and $p>2$, then $[x,y^p]=1$. Hence $[\Phi(G),x]=1$
and $K$ is abelian. If $[c,x]=a=1$ and $p=2$, then $[x,y^2]=b\neq
1$. By Lemma \ref{A_2-property} (2), $\lg x,y^2\rg\in\mathcal{A}_1$.
By \ref{G'<N}, $G'\le\lg x,y^2\rg$. Hence $K=\lg
x,y^2\rg\in\mathcal{A}_1$.

By the above argument, all maximal subgroup of $G$ are abelian or
minimal non-abelian. Hence $G\in\mathcal{A}_2$, a contradiction.
\qed

\begin{cor}
\label{cor=c=3 exp(G')=p} Suppose that $G$ is a finite
metahamiltonian $p$-group having an elementary abelian derived
group. If $c(G)=3$, then $d(G)=2$ and $p$ is odd.
\end{cor}
\demo By Theorem \ref{thm=c=3 exp(G')=p}, $G\in\mathcal{A}_2$. Then
the results follow from Corollary \ref{A_2-property} (3). \qed

\bigskip
\begin{center}
{\bf Acknowledgments}
\end{center}

\medskip

We cordially thank Professor Joseph Brennan for his assistance on the exposition and language of the paper.

\end{document}